\documentstyle[12pt]{article}

\def\R{{\bf R}}
\def\C{{\bf C}}
\def\Z{{\bf Z}}
\def\Subset{\subset\subset}
\def\supp{{\rm supp\,}}
\def\tr{{\rm tr\,}}
\def\mm{{\cal M}}
\def\ddcl{(dd^c\log|f|^2)^l}
\def\lddcl{\log|f|^2(dd^c\log|f|^2)^l}

\def\ddcq{(dd^c\log|f|^2)^q}
\def\e{\varepsilon}
\def\be{\begin{eqnarray*}}
\def\ee{\end{eqnarray*}}
\def\ben{\begin{eqnarray}}
\def\een{\end{eqnarray}}
\def\D{{\cal D}}
\def\lj{\lambda^{(j)}}
\def\<{\langle}
\def\>{\rangle}
\begin{document}

\title{Almost periodic currents, chains and divisors in tube domains}
\author{S.Yu. Favorov, A.Yu. Rashkovskii and L.I.Ronkin}
\date{}

\maketitle

Theory of almost periodic functions in one real or complex
variable was rapidly developing in the twenties-forties. The interest
to the theory was provided by both its exceptional beauty and the
significant role of almost periodic functions (uniform limits of
exponential sums with imaginary exponents) in various branches of
mathematics. The basis of the theory, its developments and
applications were founded by Bohr, Wiener, Weyl, Jessen,
Besikovich, Stepanov, Levitan, and numerous other mathematicians.
Further developments of the theory were directed mainly to its
extensions to abstract algebraic structures and applications to
differential equations. At that
time the natural problem of extension of theory of holomorphic
almost periodic functions to several complex variables faced
undeveloped proper multidimensional technique. The perspectives have
been opened by theory of closed positive currents and
Monge-Amp\`ere operators.  Recently, a series of results were obtained
on holomorphic almost periodic functions and mappings in tube
domains, and in particular on their zero set distribution (see
\cite{keyR4,keyR2,keyR3,keyR5,keyRa1,keyRa2,keyF,keyR1,keyRaRF}). 
A crucial point of the approach was an application of
distribution theory machinery that had never been applied to
holomorphic almost periodic functions before. Developing this method
it is natural to consider more general almost periodic objects such
as divisors and holomorphic chains, almost periodic distributions and
currents.  In \cite{keyR1} a notion of almost periodic distribution in a tube
domain was introduced. A definition of almost periodic divisor based
on that notion was given, and existence of density of such a divisor
was proved. It was also shown that the divisor of an almost periodic
holomorphic function is almost periodic.  The converse
relation is not true even in case of one variable \cite{keyT2}.

In the present paper a notion of almost periodic current is introduced,
as well as a notion of almost periodic holomorphic chain proceeded
from that definition. Such a chain can be defined either as a special
case of almost periodic currents or as a holomorphic chain whose
trace measure is an almost periodic distribution (as was done 
in \cite{keyR1}
for almost periodic divisors).  We prove that the both definitions
are equivalent. From the other hand, it is shown that in general
situation almost periodicity of the trace of a current does not imply
that for the current itself, even if it is closed and positive.

The zero set (regarded as a holomorphic chain) of a holomorphic
mapping can be represented as a Monge-Amp\`ere type current, and one
could expect that the zero set of an almost periodic holomorphic
mapping should be almost periodic; however we construct an example of
an almost periodic holomorphic mapping whose zero set is not almost
periodic.  Nevertheless, we prove almost periodicity of the
Monge-Amp\`ere currents corresponding to almost periodic holomorphic
mappings with certain additional properties.

Then we restrict our attention to almost periodic divisors and
construct functions that play the same role for almost periodic
divisors as the so-called Jessen functions 
(see \cite{keyJ,keyLe,keyR2,keyR3}) for almost
periodic holomorphic functions. In terms of Jessen function we give a
sufficient condition for realizability of an almost periodic divisor
as the divisor of a holomorphic almost periodic function; some
necessary condition is obtained, too.

\medskip
The paper is organized as follows. In section 1 preliminaries on
almost periodic functions and distributions are given, and almost
periodic currents are introduced and studied. Section 2 is devoted to
almost periodic holomorphic currents. Jessen functions of almost
periodic divisors are considered in section 3, and the problem of
realizability of almost periodic divisors is studied in section 4.

\section{Definition and properties of almost periodic currents}

We first recall definitions of Bohr's almost periodic function and
almost periodic distribution in a tube domain $T_G=\{z=x+iy:
x\in\R^n, y\in G\},\:G\subset\R^n$ (see\cite{keyR1}).

\medskip
{\it Definition 1.1.} A continuous function $f(z)$ in $T_G$ is called
{\it almost periodic} if for any $\e>0$ and every domain
$G'\Subset G$ the set $E_{\e,G'}(f):=\{\tau\in\R^n:
|f(z+\tau)-f(z)|<\e, \forall z\in T_{G'}\}$ is relatively dense
in $\R^n$, i.e. satisfies the condition
$$
\exists L>0,\,E_{\e,G'}(f)\cap \{t\in\R^n:|t-a|<L\}\ne\emptyset,
\quad\forall a \in \R^n.
$$

\medskip
Almost periodicity on $\R^n$ of a function $f\in C(\R^n)$ is defined
similarly by replacing the set $E_{\e,G'}(f)$ with
$$
E_{\e}(f):=\{\tau\in\R^n: |f(x+\tau)-f(x)|<\e,\quad\forall
x\in\R^n\}.
$$
To introduce almost periodic distributions we need the following
notation.

The space of functions of smoothness $\ge p,\: 0\le p\le\infty$, with
compact support in a domain $\Omega$ is denoted by ${\cal D}(\Omega,
p)$. The corresponding space of distributions (continuous linear
functionals) is denoted by ${\cal D}'(\Omega,p)$. When $p=\infty$ we
write as usual ${\cal D}(\Omega)$ and ${\cal D}'(\Omega)$ for
${\cal D}(\Omega,p)$ and ${\cal D}'(\Omega,p)$, respectively. The action
of $f\in{\cal D}'(\Omega,p)$ on a test function $\varphi\in
{\cal D}(\Omega,p)$ is denoted in a standard way by $(f,\varphi)$.

\medskip
{\it Definition 1.2.} A distribution $f\in{\cal D}'(\Omega,p),
\ p<\infty$, is called {\it almost periodic of order} $\le p$ (or just of order
$p$) if for any function $\varphi\in{\cal D}(\Omega,p)$, the function
$\left(f(z),\varphi(z-t)\right)$ is almost periodic on $\R^n$.

\medskip
A notion of almost periodic current is defined in a similar way. Let
$\Phi$ be an outer differential form on $\Omega\subset\C^n$ of
bidimension $(k,l)$, smoothness $p\ge 0$, and
$\supp\Phi\subset\Omega$. The space of all such forms is denoted by
${\cal D}_{k,l}(\Omega,p)$. The corresponding space of currents is
denoted by ${\cal D}'_{n-k,n-l}(\Omega,p)$. The action of a current $F$
on a form $\Phi$ is designated by $(F,\Phi)$.

\medskip
{\it Definition 1.3.} A current $F\in {\cal D}'_{n-k,n-l}(T_G,p)$
is called {\it almost periodic of order} $\le p$
if for any form $\Phi\in {\cal D}_{n-k,n-l}(T_G,p)$, the function
$\left(F(z),\Phi(z-t)\right)$ is almost periodic on $\R^n$.

\medskip
As is known, every current $F$ can be considered as an outer
differential form whose coefficients $F_{I,J}$ are distributions.
Comparing Definitions 1.2 and 1.3 one can see that a current $F$ is
almost periodic in $T_G$ if and only if its coefficients are almost
periodic distributions in $T_G$.

It was proved in \cite{keyR1} that almost periodic distributions in $T_G$
have properties similar to that of usual almost periodic functions.
In particular, notions of the spectrum and the Fourier series is
correctly defined, and uniqueness and approximation theorems hold
true. All this carried over  to almost periodic currents  trivially
due to the mentioned connection between almost periodicity of a
current and of its coefficients. We state here only the mean value
theorems and a theorem on compactness of translations which we will
need in future.

Usual almost periodic functions and almost periodic distributions can
be  described as functions whose translations form a compact family
in the corresponding topology. For usual almost periodic functions,
it is Bochner's definition equivalent to the original Bohr's
definition; the case of almost periodic distributions is treated in
\cite{keyR1}. A similar fact is true for currents.

\medskip
{\bf Theorem 1.1.} {\sl In order that a current
$F\in {\cal D}'_{k,l}(T_G,p)$ be almost periodic, it is necessary and
suffucient that for any sequence of the currents $F(z+h_j),\
h_j\in\R^n, \, j=1,2,\ldots$, one can choose a subsequence
$F(z+h_{j_k})$ converging in ${\cal D}'_{k,l}(T_G,p)$ uniformly on
every collection of forms $\{\Phi(z+t): t\in\R^n\},\: \Phi\in {\cal
D}_{n-k,n-l}(T_G,p)$.}

\medskip
To formulate the mean value theorems for currents, we introduce some
more notation. Let
\be
\|x\|&=&\|(x_1,\ldots,x_n)\|=\max_{1\le i\le n}|x_i|;\\
I^{(m)}=I&=&(i_1,\ldots,i_m),\,1\le i_1\le n,\ldots,1\le i_m\le n;\\
J^{(m)}=J&=&(j_1,\ldots,j_m),\,1\le j_1\le n,\ldots,1\le j_m\le n;\\
dz^I&=&dz_{i_1}\wedge \cdots \wedge dz_{i_m}, \\
d\bar z^J&=&d\bar z_{j_1}\wedge \cdots \wedge d\bar z_{j_m},\\
dx^I&=&dx_{i_1}\wedge \cdots \wedge dx_{i_m},
\ee
and let $\bar I$ denote the sample from the numbers $1,2,\ldots,n$,
complementary to $I$. Given $F\in {\cal D}'_{k,l}(\Omega,p)$,
$F_{I,J}$ denotes its coefficients considered as distributions from
${\cal D}'(\Omega,p)$, so that the current $F$ has  the
representation
$$
F=(i/2)^m\sum_{I,J}F_{I,J}dz^I\wedge d\bar z^J.
$$

We do it in the case $F\in D'_{k,l}$ as well.

The Lebesgue measure in $\R^n$ is denoted by
$m_n$, and $d^c$ denotes the operator ${1\over4i}(\partial -
\bar\partial)$.  Thus, $dd^c={i\over 2}\partial\bar\partial$.

\medskip
{\bf Theorem 1.2.} {\sl Let $F\in {\cal D}'_{m,m}(T_G,p)$ be an
almost periodic current. Then there exists a current $\mm=\mm(F)\in
{\cal D}'_{m,m}(T_G,p)$ with the coefficients
$\mm_{I,J}=e_{I,J}\otimes m_n,\: e_{I,J}\in{\cal D}'(G,p)$, such that
\ben
\lim_{\nu \to \infty}\left(\frac1{2\nu}\right)^n\int_{\|t\|<\nu}
(F(z+t),\Phi(z))dx=(\mm,\Phi),\\ 
\forall\Phi\in{cal D}_{n-m,n-m}(T_G).\nonumber
\een
If $p=0$ and the coefficients of a form $\Phi$ are functions from
${\cal D}(G,0)$, then in addition to (1), the equality holds}
$$
\exists\,\lim_{\nu\to\infty}\left(\frac1{2\nu}\right)^n\int_{\|x\|<\nu,y\in G}
F(z)\wedge\Phi(z)=(\mm,\Phi).
$$

\medskip
Another variant of the mean value theorem is given by

\medskip
{\bf Theorem 1.3.} {\sl Let $F$
be an almost periodic current from $ {\cal D}'_{m,m}(T_G,0)$ and 
$$
F^{(\nu)}=(i/2)^{n-m}\sum_{I,J}F_{I,J}(\nu x+iy)dz^I\wedge d\bar z^J.
$$
Then the equality
$$
\lim_{\nu\to\infty}F^{(\nu)}=\mm (F)
$$
takes place in $ {\cal D}'_{m,m}(T_G,0)$. }

\medskip
Note also that the operator $F\mapsto\mm(F)$ keeps such properties of
currents as positivity and closedness.

As was said above, Theorems 1.2 and 1.3 are trivial consequences of
the corresponding results from \cite{keyR1} on almost periodic distributions.
Non-trivial points appear when applying almost periodic currents to
the study of almost periodic holomorphic mappings, divisors, and
holomorphic chains. Such problems will be treated in the following
sections.

\section{Almost periodic holomorphic mappings and chains}

{\bf 1.} Important examples of almost periodic currents appear from
almost periodic holomorphic mappings.

Let $f:T_G\to\C^q, q\le n$, be an almost periodic holomorphic
mapping,  that is $f=(f_1,\ldots,f_q)$ where $f_1,\ldots,f_q$ are
almost periodic holomorphic functions in $T_G$. It is natural to
expect that its zero set - or, more precisely, the holomorphic chain
$Z_f$ generated by the mapping - inherits the property of almost
periodicity. Actually the things are not so good, and we will show
below that the zero set of an almost periodic holomorphic mapping
need not be  almost periodic. From the other hand, we will prove that
under certain natural conditions on a mapping, the corresponding
current of integration over $Z_f$ has to be almost periodic, as well
as some other currents generated by $f$. Then we come to the notion
of almost periodic holomorphic current and establish some its
properties.

Recall some notions and facts concerning objects under consideration.

\medskip
{\it Definition 2.1.} A {\it holomorphic chain} of dimension $m$ (or
codimension $n-m$), $0\le m\le n$, in a domain $\Omega\in\C^n$ is a
pair $Z=(|Z|,\gamma_Z)$, where $|Z|$ is an analytic set in $\Omega$ of
pure dimension $m$, and $\gamma_Z$ is a function on $|Z|$ that has
constant integer value on each connected component of the set
$reg\,|Z|$ of regular points of the set $|Z|$. If $|Z|=f^{-1}(0)$ for
a holomorphic mapping $f$, and $\gamma_Z(z)$ is the
multiplicity of the mapping $f$ in the point $z$, then the chain is
said to be given by the mapping $f$ and is denoted by $Z_f$, as well
as $\gamma_Z(z)$ is replaced with $\gamma_f(z)$. In  case of $\dim
Z=n-1$ the holomorphic chain is called {\it divisor}.

\medskip
It is clear that a holomorphic chain can also be considered as a
mapping from $\Omega$ to $\bf Z$ with the specified properties of its
support and values.

In the sequel we assume without special mentioning that $\gamma_Z>0$.
The current of integration over the chain $Z$ is denoted by $[Z]$.
Such a current is defined by the equation
$$
([Z],\Phi)=\int_{reg|Z|} \gamma_Z(z)\Phi\mid_{|Z|}, \quad 
\Phi\in{\cal D}_{m,m}(\Omega,0),
$$
where $\Phi\mid_{|Z|}$ is the restriction of the form $\Phi$ to
$reg\,|Z|$.

The space of holomorphic functions on a domain $\Omega$ is denoted as
usual by $H(\Omega)$.

Let a holomorphic mapping $f:\Omega\to\C^q,\;q\le n$, satisfy the
condition $\dim |Z_f|\le n-q$ (that is either $|Z_f|=\emptyset$ or
$\dim_a |Z_f|=n-q,\;\forall a\in |Z_f|$). Then as is known (see
\cite{keyGrK}), the Monge-Amp\`ere type currents $\ddcl$ 
and $\lddcl,\; l<q$,
have locally summable coefficients, and the current
$$
\ddcq:=dd^c(\log|f|^2(dd^c\log|f|^2)^{q-1})
$$
coincides (up to a constant factor) with the current $[Z_f]$. The
following statement was proved in \cite{keyR2} for $l=q-1$ and in 
\cite{keyRa1}
for $l<q$: if a sequence of holomorphic mappings $f_j:\Omega\to\C^q$
converges to the mapping $f$ as $j\to\infty$, uniformly on every
compact subset of $\Omega$, then
$$
\lim_{j\to\infty}\log|f_j|^2(dd^c\log|f_j|^2)^l=\lddcl,\quad\forall l<q
$$
in the space ${\cal D}_{l,l}$.

Let $f$ be an almost periodic holomorphic mapping of a domain $T_G$
into $\C^q,\;q\le n$. Then according to one of the equivalent
definitions of almost periodic function (Bochner's definition), for
any sequence $\{h_j\}\subset\R^n$ one can choose a subsequence
$\{h_{j_k}\}$ such that  $f(z+h_{j_k})$ converge as $k\to\infty$,
uniformly on every tube domain $T_{G'},\;G'\Subset G$.  The
collection of all such limit mappings is denoted by ${\cal F}_f$. The
mapping $f$ is called {\it regular} (see \cite{keyR2,keyR3}) 
if $\dim |Z_{\tilde f}|\le n-q,\;\forall\tilde f\in{\cal F}_f$.  
The set of all regular almost periodic holomorphic mappings of the 
domain $T_G$ into $\C^q$ is denoted by $R_q(G)$. A sufficient condition 
for an almost periodic holomorphic mapping to be regular, in terms of 
its spectrum, is given in \cite{keyR2} (see \cite{keyR3} also).

Now we have all necessary to formulate the mentioned statement on
almost periodic currents connected with almost periodic holomorphic mappings.

\medskip
{\bf Theorem 2.1.} {\sl If $f\in R_q(G)$ then the currents $\lddcl,\quad
l<q$, and $\ddcl,\ l\le q$, are almost periodic in $T_G$.}

\medskip
{\it Proof.} Differentiation evidently keeps
almost periodicity of a current. So the equality
$$
\ddcl=dd^c(\log|f|^2(dd^c\log|f|^2)^{l-1})
$$
shows that it suffices to prove the statement of the theorem for
the currents
$$
A_f^{(l)}:=\lddcl, \quad l<q.
$$

Let $\Phi\in{\cal D}_{n-l,n-l}(T_G,0)$. Denote
$Q(t)=\left(A_f^{(l)}(z+t), \Phi(z)\right)$, and let $\{h_j\}$ be an
arbitrary sequence in $\R^n$. By almost periodicity of the mapping
$f$, one may take without restricting any generality $f(z+h_j)$
converging as $j\to\infty$ to some mapping $\tilde f\in{\cal F}_f$,
uniformly on every domain $T_{G'},\; G'\Subset G$. Since ${\tilde f}\in
R_q(G)$, $\dim |Z_{\tilde f}|\le n-q$ and thus the current $A_f^{(l)}$ is well
defined. We claim that for these $h_j$,
$$
Q(t+h_j)\to \tilde Q(t):=(A_{\tilde f}^{(l)}(z+t),\Phi(z)).
$$
uniformly in $\R^n$.

Supposing the contrary, there exist a number $C_0$ and a subsequence
$\{t_j\}\subset\R^n$ such that $|Q(t_j+h_j)-\tilde Q(t_j)|\ge C_0,\;
\forall j$. Passing if necesary to a subsequence one may take that,
uniformly in $T_{G'},\; \forall G'\Subset G$, there exists the limit
$$
\lim_{j\to\infty}{\tilde f}(z+t_j)=:\hat f(z).
$$
It is clear that then $f(z+t_j+h_j)\to\hat f(z)$ and $\hat f\in
R_q(G)$. By the mentioned convergence theorem for the Monge-Amp\`ere
type currents, $A_{\tilde f(z+t_j)}^{(l)} \to A_{\hat f}^{(l)}(z)$
and $A_{f(z+h_j+t_j)}^{(l)} \to A_{\hat f}^{(l)}(z)$.
Therefore
\be
\limsup_{j\to\infty}|Q(t_j+h_j)-\tilde Q(t_j)|\le
\limsup_{j\to\infty}|(A^{(l)}_{f(z+t_j+h_j)},\Phi)-
(A^{(l)}_{\tilde f},\Phi)|+\\
\limsup_{j\to\infty}|(A^{(l)}_{\tilde f(z+t_j)},\Phi)-
(A^{(l)}_{\hat f},\Phi)|=0.
\ee
that contradicts to the choice of the sequence $\{t_j\}$. So, our
claim is proved, and then the function $Q(t)$ is almost periodic on
$\R^n$ and it is true for every form $\Phi\in{\cal
D}_{n-l,n-l}(T_G,0)$. Thus the current $A_f^{(l)}$ is almost
periodic. The theorem is proved.

\medskip
The following two results are easy consequences of Theorem 2.1 and
the foregoing exposition.

\medskip
{\bf Corollary 2.1.} {\sl The integration current $[Z_f]$ over the
holomorphic chain generated by a mapping $f\in R_q(G)$, is almost periodic.}

\medskip
{\bf Corollary 2.2.} {\sl If $f\in R_q(G)$ then the currents
$A_f^{(l-1)}$ and $a_f^{(l)}:=\ddcl,\: l\le q$, have mean values
$$
\mm(A^{(l-1)}_f)=:\tilde A^{(l-1)}_f
$$
and}
\ben
\mm(a^{(l)}_f)=:\tilde a^{(l)}_f.
\een

\medskip

The currents $\tilde a_f^{(l)}$ inheriting properties of the current
$a_f^{(l)}$ , are evidently positive, closed, and connected with the
currents $\tilde A_f^{(l-1)}$ by the relation
\ben
dd^c\tilde A^{(l)}_f=\tilde a^{(l)}_f.
\een
Corollary 2.2 was proved earlier as an original result not exploiting
almost periodicity of the corresponding currents, in \cite{keyR2} for $l=q-1$
and in \cite{keyRa1} for any $l$.

Now we construct an example of almost periodic holomorphic mapping
showing that the {\it condition} $f\in R_q(G)$ {\it is essential} for these
results.

For $k\in{\bf N}$ set
$$
g_k(\zeta)={\sin\pi\zeta\over\sin\frac{\pi\zeta}{k}}.
$$
It is an entire periodic function in the complex plane, and for an
integer $n$, $g_k(n)=0$ if and only if $k$ is not a divisor of $|n|$
and $n\neq 0$. Choose a sequence $\{a_k\}\subset\R$ such that
$$
0<a_k\sup\{|g_k(\zeta)|: |{\rm Im}\,\zeta|<3\}\cdot
\sup\{|\sin\pi\zeta)|: |{\rm Im}\,\zeta|<3\}  <k^{-2}
$$
and consider the mapping $f(z_1,z_2)=(f_1,f_2)$:
\be
f_1(z_1,z_2) & = & \sin {\pi(z_1-2)\over 5},\cr
f_2(z_1,z_2) & = & \sum_{k\ge 2} a_kg_k(z_1)\sin\pi k z_2.
\ee
By the choice of $a_k$, the function $f_2$ is defined by the series
that converges uniformly in the tube domain $T_G=\{z\in\C^2:\:
|{\rm Im}\,z_j|<2,\: j=1,2\}$ and thus is an almost periodic
holomorphic  function there. So, $f$ is an almost periodic
holomorphic mapping from $T_G$ to $\C^2$.

First we show that $|Z_f|$ is discrete. To this end it suffices to
prove that $|Z_f|$ contains no set of the form $\{z:\: z_1=5n+2\},\;
n\in{\bf Z}$. Let $m>1$ be an integer, then $g_k(m)=0 \ \forall k>m$,
and so
$$
f_2(m,z_2)  =  \sum_{k\in A(m)} a_kg_k(m)\sin\pi k z_2,
$$
where $A(m)$ is the set of all divisors $k>1$ of the number $m$.
Since $g_k(m)\neq 0\ \forall k\in A(m)$, $f_2(m,z_2)\not\equiv 0$ for
any $m>1$. It remains true for $m<-1$ as $f_2(-m,z_2) =f_2(m,z_2)$.
It proves our claim because $|5n+2|>1,\:\forall n\in{\bf Z}$.

Now we show that the current $[Z_f]$ is not almost periodic. Let $p>1$ be a
prime number,  then $A(p)=\{p\}$ and
$$
f_2(p,z_2)  =  a_pg_p(p)\sin\pi p z_2.
$$
Therefore $f_2$ vanishes in the points $(p,l/p),\:
l=0,\pm 1,\ldots,\pm p$. Let $P$ denote the set of all primes of
the form $5n+2$. Then we have
$$
\{(p,l/p):\: p\in P,\: l=0,\pm 1,\ldots,\pm p\}\subset Z_f.
$$
Hence any open ball of radius one and center at $(p,0),\: p\in P$,
contains at least $2p-1$ points of the set $Z_f$. Since every such
ball is a translation of the ball $B=\{z:\: |z|<1\}$ and the set
$P$ is infinite, the function $\left([Z_f],\varphi (z-t)\right)$ is
not bounded on $\R^2$ for a non-negative test function $\varphi,\;
{\rm supp}\,\varphi\supset B$, so that the current $[Z_f]$ is not
almost periodic.

\bigskip
2. Now we introduce the following

{\it Definition 2.2.} A holomorphic chain $Z$ in a domain $T_G$ is
called almost periodic if the current $[Z]$ of integration over $Z$ is
almost periodic\footnote{In \cite{keyR1} for the case of ${\rm codim}\,Z=1$,
i.e.  when $Z$ is a divisor, a definition of almost periodicity was
given, different from Definition 2.2 and based on the notion of
almost periodic distribution.  For the relation between the two
definitions, see Theorem 2.3 below.  In the case $n=1$ a definition
of almost periodic divisor (almost periodic set) was given in
\cite{keyKL,keyLe} under certain restrictions on divisors, 
and in \cite{keyT1} in general
situation; that definition is based on an approach different from
considered here.}.

\medskip
Corollary 2.2 implies that the chain $Z_f$ associated to a mapping
$f\in R_q(G)$ is almost periodic.

By the mentioned properties of almost periodic currents, the
integration current over an almost periodic holomorphic chain has
the mean value. Now we study the geometrical sense of this mean
value.

Let $Z$ be an almost periodic holomorphic chain of dimension $l$
in $T_G$. Denote $\beta_l= {1\over l!}(dd^c |z|^2)^l$, and let
$\chi_E=\chi_E(z)$ be the indicator function of a Borelian set $E$.
Then (see for example \cite{keyLG}) the value
\be
V_Z(\Omega):=(\chi_{\Omega}[Z],\beta_l)=
\int\chi_{\Omega}[Z]\wedge\beta_l=\int_{|Z|}\gamma_Z(z)\beta_l,
\ee
where $\Omega\Subset T_G$, is equal to the $2l$-dimensional volume of
the chain $Z$ in the domain $\Omega$. By ${\rm tr}\, [Z]$ denote the trace
measure of the current $[Z]$; in other words,
$$
\tr[Z]=\sum_I[Z]_{I,I}.  
$$
Evidently,
$$
[Z]\wedge\beta_l=\tr[Z]\beta_n.
$$
According to Theorem 1.3 we construct the currents $[Z]^{(\nu)}$ and
then the measures
$$
\tr[Z]^{(\nu)}=:\mu_Z^{(\nu)}.
$$
Since the current $[Z]$ is positive, the measures $\mu_Z^{(\nu)}$ are
positive, too. Denote
$$
\Pi_{t,E}=\{z:\|x\|<t,y\in G'\}, \quad  E\subset G,
$$
 and observe that for $G'\Subset G$
\ben
\mu_Z^{(\nu)}(\Pi_{1,G'})=\int_{\Pi_{1,G'}}\tr[Z]^{(\nu)}\beta_n=
\nu^{-n}\int_{\Pi_{\nu,G'}}[Z]\wedge\beta_l=\nu^{-n}V_Z(\Pi_{\nu,G'}).
\een
By Theorem 1.3, the currents $[Z]^{(\nu)}$ converge in ${\cal
D}'_{l,l}(T_G,0)$ to a current $\mm([Z])$  of the form
$$
\mm([Z])=(i/2)^{n-l}\sum_{I,J}b_{I,J}\otimes dz^I\wedge d\bar z^J,
$$
where $I=I^{(n-l)},\: J=J^{(n-l)},\: b_{I.J}\in {\cal
D}'(G,0)$, and $m_n$ is the Lebesgue measure in $\R_{(y)}^n$.

Define a measure $\mu_Z$ on $G$ by setting for Borelian sets
$E\subset G$
$$
\mu_Z(E)=\int_{\Pi_{1,E}}\tr\mm([Z])\beta_n.
$$
so that
$$
\mu_Z=\sum_Ib_{I,I}
$$
and $\mu_Z\otimes m_n={\rm tr}\, \mm([Z])$.

The convergence of the currents $[Z]^{(\nu)}$ to the current $\mm([Z])$
implies that the measures $\mu_Z^{(\nu)}$ converge weakly to the
measure $\mu_Z\otimes m_n$ as $\nu\to\infty$. Therefore (see for
example \cite{keyLa,keyR3}) if a domain $G'\Subset G$ satisfies
the conditions
$$
(\mu_Z\otimes m_n)(\partial\Pi_{1,G'})=0
$$
or, that is the same, $\mu_Z(\partial G')=0$, then
$$
\lim_{\nu\to\infty}\mu_Z^{(\nu)}(\Pi_{1,G'})=
(\mu_Z\otimes m_n)(\Pi_{1,G'})=2^n\mu_Z(G').
$$
Combined with (4) it gives us the following statement.

\medskip
{\bf Theorem 2.2.} {\sl Let $Z$ be an almost periodic holomorphic
chain of a dimension $l<n$ in a domain $T_G$ and let a domain
$G'\Subset G$ satisfy $\mu_Z(\partial G')=0$. Then there exists the
limit $\lim_{\nu\to\infty}(2\nu)^{-n} V_Z(\Pi_{\nu,G'})$ and the
equation takes place}
\ben
\lim_{\nu\to\infty}(2\nu)^{-n} V_Z(\Pi_{\nu,G'})=\mu_Z(G').
\een

\medskip
It is natural to call the measure $\mu_Z$ from Theorem 2.2 {\it the
density of the chain} $Z$, so that the theorem can be viewed as a
theorem on existence and evaluation of the density of an almost
periodic chain. In the case of $Z=Z_f,\: f\in R_q(G)$, the result
was obtained earlier (with no regard to almost periodicity of the
current $[Z_f]$) in \cite{keyR4} for $q=1$ and in \cite{keyR2} for $1\le 
q\le n$ (see \cite{keyR3} also).
Those results were preceded by the classical Jessen theorem (\cite{keyJ},
see also \cite{keyLe,keyLt}) on density of the zero set of an almost periodic
holomorphic function of one complex variable.

The following useful property of almost periodic holomorphic chains
is connected with their densities.

\medskip
{\bf Theorem 2.3.} {\sl Let $Z$ be an almost periodic holomorphic
chain in $T_G$ and let $\mu_Z(G')=0$ for some domain $G'\subset G$.
Then $|Z|\cap T_{G'}= \emptyset$.}

\medskip
{\it Proof.} Suppose the contrary. Then there exists a form
$\Phi\in {\cal D}_{l,l}(T_G,0)$ with $l=\dim Z$ such that
$([Z],\Phi)\neq 0$. It follows from almost periodicity of the
function $\left([Z],\Phi(z-t)\right)$ that for a certain $L$ the
function is identically zero on no domain $\{t\in\R^n:\:
||t-kL||<L\}, \: k\in {\bf Z}^n$. Therefore if a domain $G"$
satisfies the condition ${\rm supp}\,\Phi\subset G"\Subset G'$ then
for any $k\in {\bf Z}^n$ there is a point $z^{(k)}\in |Z|\cap
\{z\in\C^n:\: ||x-kL||<L,\, y\in G"\}$. Choose $\e={\rm dist}\,
(G", \partial G')$ and set $\omega_\e(z^{(k)})=\{z\in\C:\:
|z-z^{(k)}|<\e\}$. By the Lelong theorem on the lower bound for
the mass of a closed positive current (see \cite{keyLG}),
$$
V_Z(\omega_{\e}(z^{(k)}))\ge{\rm const\,}\e^{2l},
$$
and thus for $\nu$ sufficiently great,
$$
V_Z(\Pi_{\nu,G'})\ge{\rm const\,}\e^{2l}(2\nu)^l.
$$
In view of (5) it contradicts to the assertion $\mu_Z(G')=0$. The
theorem is proved.

\bigskip
3. As was mentioned, almost periodicity of a current $F$ is
equivalent to that of its coefficients $F_{I,J}$. In case of $F=[Z]$,
$Z$ being a holomorphic chain, a weaker assumption is sufficient for
almost periodicity of the current (and so for the chain).

\medskip
{\bf Theorem 2.4.} {\sl Let $Z$ be a holomorphic chain of dimension
$l$ in $T_G$, and let the distribution ${\rm tr}\, [Z]$ is almost periodic
in $T_G$. Then the chain $Z$ is almost periodic, too.}

\medskip
This theorem allows us to define almost periodic holomorphic chain
as  such a chain that the corresponding integration current has
almost periodic trace. Actually, for $l=n-1$ it is the mentioned
definition of almost periodic divisor from \cite{keyR1} based on
consideration of the trace of the current of integration over the
divisor.

\medskip
{\it Proof.} Denote by $[Z]_t$ the "translated" current, i.~e. the
current acting on a form $\Phi\in{\cal D}_{l,l}(T_G,0)$ as
$([Z]_t,\Phi)=([Z],\Phi(z-t))$. The theorem would be proved if we
show that for any sequence $\{t_j\}\subset\R^n$ one can choose a
subsequence $\{t_{j_k}\}$ such that  $\{[Z]_{t_{j_k}+t}\}$
converges as $k\to\infty$, uniformly on $\R^n,\: \forall
\Phi\in{\cal D}_{l,l}(T_G,0)$.

Note that as follows from the properties of almost periodic
distributions (see \cite{keyR1}), the value
$$
\|Z\|_{G^0}:=\sup_{a\in\R^n}\int_{a+\Pi_{1,G^0}}\tr [Z]\beta_n
$$
is finite for every domain
$G^0\Subset G$. Therefore the masses of the "trasnslated" currents
are uniformly bounded on each compact subset $K$ of $T_G$. Hence the
family $\{[Z]_t\}$ is weakly compact and thus, given a sequence
$\{t_j\}\subset\R^n$, one can choose a subsequence $\{t_{j'}\}$ such
that the currents $[Z]_{t_{j'}}$ form a converging sequence in
${\cal D}'_{l,l}(T_G,0)$. At the same time the assumed almost
periodicity of ${\rm tr}\,[Z]$ provides the existence of a
subsequence $\{t_{j"}\}$ of $\{t_{j'}\}$ such that the functions
$({\rm tr}\,[Z]_{t+t_{j"}},\varphi)$ or, that is the same, the
functions $([Z]_{t+t_{j"}},\varphi\beta_l)$, converge as
$j"\to\infty$, uniformly on $\R^n,\:\forall\varphi\in
{\cal D}(T_G,0)$. For the sake of brevity we will write $t_j$ for
$t_{j"}$. So, we deal with the sequence $\{t_j\}\subset\R^n$
sarisfying the conditions
\ben
\exists \lim_{j\to\infty}[Z]_{t_j}=:F ;\nonumber\\
\ \\
\lim_{j\to\infty} \sup_t |([Z]_{t+t_{j"}}-F_t,\varphi\beta_l)|=0,
\quad \forall \varphi\in {\cal D}(T_G,0);\nonumber
\een
here $F_t=F(z+t)$.

We claim that these conditions imply the equality
$$
\lim_{j\to\infty} \sup_{t\in\R^n}
|([Z]_{t+t_{j"}}-F_t,\Phi)| =0,
\quad \forall \Phi\in {\cal D}_{l,l}(T_G,0).
$$
Suppose the contrary. Then for some $\e>0$ and some form
$\Phi\in {\cal D}_{l,l}(T_G,0)$ there exists a sequence $\{h_j\}\in
\R^n$ such that
\ben
|([Z]_{t_j+h_j}-F_{h_j},\Phi)|>\e, \quad \forall j.
\een
Choose a subsequence $\{j'\}\subset\{j\}$ in a way the equalities
\ben
\exists \lim_{j'\to\infty}[Z]_{t_{j'}+h{j'}}:=\tilde F\nonumber\\
\ \\
\exists \lim_{j'\to\infty}F_{h_{j'}}=\hat F\nonumber
\een
hold. It is possible due to the mentioned weak compactness of the
family $\{[Z]_t\}_{t\in\R^n}$ and so of the family
$\{F_t\}_{t\in\R^n}$.

Passing to the limit in (7)  as $j'\to\infty$ we get
\ben
|(\tilde F-\hat F,\Phi)|>\e.
\een
However it is impossible. Indeed,
\be
|([Z]_{t_{j'}+h_{j'}},\beta_l\varphi)-(F_{h_{j'}},\beta_l\varphi)|=
|([Z]_{t_{j'}}-F,\beta_l(z)\varphi(z-h_{j'}))|\le\\
\sup_{t\in\R^n}|([Z]_{t_j}-F,\ \beta_l(z)\varphi(z-t))| 
\ee
for all $j$, for all $\varphi\in{\cal D}(T_G,0)$. Together with (6) 
it implies that for $j'\to\infty$
$$
|([Z]_{t_{j'}+h_{j'}}\wedge\beta_{n-l},\varphi)-
(F_{h_{j'}}\wedge\beta_{n-l},\varphi)|\to 0.
$$
In its turn it gives us by (8)
$$
(\tilde F\wedge\beta_{n-l},\varphi)=(\hat F\wedge\beta_{n-l},\varphi) \quad 
\forall \varphi\in\D(T_G,0).
$$
Therefore
\ben
{\rm tr}\hat F={\rm tr}\tilde F.
\een
Each of the currents $[Z]_t$ is the integration current over the
corresponding holomorphic chain. Since their masses are uniformly
bounded, their limits are also currents of integration over
holomorphic chains due to the known Federer-Bishop theorem (see for
example \cite{keyH}). So, there exist holomorphic chains $\tilde Z$ and
$\hat Z$ of dimension $l$ such that $\tilde F=[\tilde Z]$ and $\hat
F=[\hat Z]$. By (10), these chains satisfy the condition ${\rm tr}\,
[\tilde Z]={\rm tr}\,[\hat Z]$. However for integration currents over
holomorphic chains the coincidence of the traces implies coincidence
of the currents themselves. Therefore $(\tilde F-\hat F,\Phi)=0 \
\forall\Phi\in{\cal D}_{l,l}(T_G,0)$ that contradicts to (9). The
theorem is proved.

\medskip
In conclusion of the section we give an example showing that the
trace of a closed positive current can be almost periodic while the current
itself need not be almost periodic.

Set $u(z)=u(z_1,z_2)=y_1^2+y_2^2+x_1{\rm Re}\,e^{z_2^2}$. Then
$$
\left(\frac{\partial^2 u}{\partial z_i\partial\bar z_j}\right)=
\pmatrix{2 & -2\bar z_2 e^{\bar z_2^2} \cr -2z_2 e^{z_2^2} & 2 \cr}
$$
and
$$
\det \left(\frac{\partial^2 u}{\partial z_i\partial\bar z_j}\right)
=4\left(1-(x_2^2+y_2^2)e^{-2(x_2^2-y_2^2)}\right).
$$
Elementary calculations with a fixed $y_2\in\left( -{1\over\sqrt
2}, {1\over\sqrt 2}\right)$ show that
$$
\sup_{-\infty<x_2<\infty} (x_2^2+y_2^2)e^{-2(x_2^2-y_2^2)} =
{1\over 2}e^{4y_2^2-1}.
$$
So for a fixed $y_2\in\left( -{1\over 2}, {1\over 2}\right)$,
$$
\det \left(\frac{\partial^2 u}{\partial z_i\partial\bar z_j}\right)>0,
$$
and thus the Levy form of the function $u$ is positively defined in
the domain $T_G\subset\C^2$ with the base $G=\{(y_1,y_2):\:
-\infty<y_1<\infty, -{1\over 2}<y_2< {1\over 2} \}$.

Let $F=(dd^c u)^2$. It is a closed positive current in $T_G$. Its
trace ${\rm tr}\, F=4m_4$ and hence is an almost periodic measure.
At the same time, if $\Phi={i\over 2}\varphi\,d z_1\wedge d\bar z_2$
with $\varphi \in{\cal D}(T_G,0)$, then the function
$(F(z+t),\Phi(z)),\:  t\in\R^2$, is not almost periodic. Indeed,
$$
(F(z+t),\Phi(z))=\int(z_2+t_2)^2e^{-(z_2+t_2)^2}\varphi(z)dm_4
$$
so that $(F(z+t),\Phi(z))\to 0$ as $t_2\to\infty$.
If the function $(F(z+t),\Phi(z))$ were almost periodic it would imply it
vanishes everywhere, that is deliberately not the case. Therefore the
current $F$ is not almost periodic.

\section{Jessen functions of almost periodic divisors}

Let $f$ be an almost periodic holomorphic function in $T_G$. By
Theorem 2.2 and Corollary 2.1, its divisor $Z_f$ has the density
$\mu_f:=\mu_{Z_f}$ connected by the mean value $\mm([Z_f])$ by the
relation
\ben
\mu_f\otimes m_n={\rm tr}\mm([Z_f])
\een
At the same time, $\mm([Z_f])=\tilde a_f^{(1)}$, so that (see (2))
\ben
\mm([Z_f])=dd^c\tilde A^{(0)}_f
\een
It was proved in \cite{keyR4}, se also \cite{keyR3}, that, given 
an almost holomorphic function $f(z)$, there exists the limit
$$
\lim_{\nu\to\infty}\left({1\over{2\nu}}\right)^n\int_{|x|<\nu}\log|f(x+iy)|dx
=:A_f(y)
$$
and the function $A_f(y)$ is convex. Following \cite{keyR4,keyR3} we will call it
{\it the Jessen function} of the function $f(z)$.\footnote{B.~Jessen
who introduced such a function in \cite{keyJ} for $n=1$, called it Jensen's
function, referring to the case where the function $f$ is periodic
and $A_f$ can be transformed into the corresponding term of the
Jensen formula.  We believe, however, that due to the deep Jessen's
investigations of specific features of the function $A_f$ in case of
his own name.}

Since $A_f(y)={1\over 2}\tilde A_f^{(0)}$, relations (11) and (12)
imply that
\ben
\mu_f=2A_f(y).
\een
so that $\mu_f$  coincides up to a constant factor with the measure
associated by Riesz (the Riesz measure) of a convex and thus
subharmonic, function $A_f(y)$. Note that, given a positive measure,
there is a family of subharmonic functions whose Riesz measures
coincide with it (each of the functions differs from other by a
harmonic function). However there may be no convex function in the
family. So, relation (13) can be considered as a certain
characteristics of the density of a divisor $Z_f$.

It was shown in the previous section that not only the divisor of an
almost periodic holomorphic function but also any almost periodic
divisor has the density. And the following question is quite natural:
do such densities possess the above property of the densities of the
divisors $A_f$? The affirmative answer is a consequence of the
following result.

{\bf Theorem 3.1.} {\sl Let $F$ be a closed positive almost periodic
current from ${\cal D}'_{(n-1,n-1)}(T_G,0)$. Then there exist a
convex function $A(y)$ in $G$ and constants $c_{j,k}\in\R,\:
j,k=1,\ldots,n$, such that}
$$
\mm(F)=dd^c(A(y)+\sum_{j<k}2c_{j,k}(x_iy_j-y_ix_j)),
$$

\medskip
{\it Proof.} Denote $\theta=\mm(F)$. As was mentioned above, the
property of a current to be closed and positive is inhereted by its
mean value. So, the current $\theta$ is closed and positive. It can
be represented in the form $\theta=\theta'+i\theta"$, where
$\theta'$ and $\theta"$ are the currents whose coefficients
$\theta'_{j,k}$ and $\theta"_{j,k}$ are real measures. Since the
current $\theta$ is positive, $\theta'_{j,k}=\theta'_{k,j},\
\theta'_{j,j}\ge 0,\ \theta"_{j,k}=-\theta"_{k,j},\
\theta"_{j,j}= 0$. Furthermore, as the measures $\theta_{j,k}=
\theta'_{j,k}+ i\theta"_{j,k}$ have the form $\theta_{j,k}=
\hat\theta_{j,k}\otimes m_n$, where $\hat\theta_{j,k}=
\hat\theta'_{j,k}+ i\hat\theta"_{j,k}$ are complex measures in $G$,
then
$$
{\partial\theta_{j,k}\over\partial x_s-}=0
$$
so that
\ben
{\partial\theta_{j,k}\over\partial z_m}=
-i/2 {\partial\theta_{j,k}\over\partial y_m},\qquad
{\partial\theta_{j,k}\over\partial \bar z_m}=
i/2{\partial\theta_{j,k}\over\partial y_m}
\een
The corresponding relations are certainly valid for the measures
$\theta'_{j,k}$ and $\theta"_{j,k}$, too. Notice further that the
currents $\theta'$ and $\theta"$ are closed. By (14), it gives us
\ben
{{\partial\theta_{j,k}'}\over{\partial y_l}}=
{{\partial\theta_{k,l}'}\over{\partial y_j}}=
{{\partial\theta_{l,j}'}\over{\partial y_k}},
\een
\ben
{{\partial\theta_{j,k}''}\over{\partial y_l}}=
{{\partial\theta_{k,l}''}\over{\partial y_j}}=
{{\partial\theta_{l,j}''}\over{\partial y_k}},\quad\forall j,k,l.
\een

Since $\theta"_{j,k}=-\theta"_{k,j}$, it follows from (16) that
$$
{{\partial\theta_{j,k}''}\over{\partial y_l}}=
-{{\partial\theta_{k,j}''}\over{\partial y_l}}=
-{{\partial\theta_{l,j}''}\over{\partial y_k}}=
-{{\partial\theta_{j,k}''}\over{\partial y_l}},
\
\forall j,k,l.
$$
Therefore
$$
{{\partial\theta_{j,k}''}\over{\partial y_l}}=0,
\quad\forall j,k,l,
$$
so that $\theta"_{j,k}=c_{j,k} m_n\otimes m_n$, where $c_{j,k}$ are
some real constants. Hence
$$
\theta''=idd^c(\sum_{j<k}c_{j,k}(z_j\bar z_k-z_k\bar z_j)).
$$

Now to prove the theorem, it is sufficient to find a convex function
$A(y)$ in $G$ such that
\ben
{{\partial^2 A}\over{\partial y_j\partial y_k}}=4\theta'_{j,k},
\quad\forall j,k.
\een
We will first solve this problem in the balls
$B(x_0,R)=\{x\in\R^n:|x-x_0|<R\}\Subset G$ and then will "glue" the solutions.
To simplify notation, we take $x_0=0$ and denote $B(0,R)$ by $B(R)$.

The desired function will be constructed by a method close to one of
the proof of Theorem 2.28 from \cite{keyLG}.

We first consider the  function
$$
A_1(y)=-\int_{B(R+\e)}|\zeta-y|^{2-n}d\sigma(\zeta),
$$
subharmonic in $\R^n$, where $\e$ is such that $B(R+2\e)
\subset G$, the measure $\sigma(\zeta)$ is defined by the equality
$$
\sigma(\zeta)=\kappa_n\sum_{l=1}^n\theta'_{l,l},
$$
and the norming constant $\kappa_n$ is determined by the condition
$$
\Delta A_1(y)=4\sum_{l=1}^n\theta'_{l,l}(y).
$$
Set
\ben
v_{j,k}(y)=\theta'{j,k}(y)-
{1\over4}{\partial^2 A_1(y)\over\partial y_j\partial y_k}.
\een
In view of (18) and (15), we have in $B(R+\e)$
$$
\sum_{l=1}^n{\partial^2v_{j,k}\over\partial y^2_l}=
\sum_{l=1}^n{\partial^2\theta'{l,l}\over\partial y_j\partial y_k}-
{1\over4}{\partial^2 \Delta A_1(y)\over\partial y_j\partial y_k}=0,
$$
that is the functions $v_{j,k}$ are harmonic in $B(R+\e)$.
Therefore, in this ball they can be represented as
$$
v_{j,k}=\sum_{s=0}^{\infty}P_{s,j,k}(y)
$$
where $P_{s,j,k}$ are homogeneous polynomials of degree $s$. It
follows directly from the definition of the functions $v_{j,k}$ that
$$
{{\partial v_{j,k}}\over{\partial y_l}}=
{{\partial v_{k,l}}\over{\partial y_j}}=
{{\partial v_{l,j}}\over{\partial y_k}},\qquad\forall j,k,l.
$$

It is easy to see that the same relations have place for the
functions $P_{s,j,k}(y)$. Together with the Euler equation
$$
\sum_{l=1}^n{{\partial P_{s,j,k}}\over{\partial y_l}}y_l=s\,P_{s,j,k}
$$
it implies that the function
$$
A_2(y)=4\sum_{j,k=1}^n\sum_{s=0}^{\infty}(s+1)^{-1}(s+2)^{-1}
P_{s,j,k}(y)y_jy_k
$$
satisfies the equations
$$
{{\partial^2 A_2(y)}\over{\partial y_j\partial y_k}}=4v_{j,k},
\quad\forall j,k,
$$
in the ball $B(R+\e)$.
Then the function $A(y):=A_1(y)+A_2(y)$ satisfies the equations
(17) in the ball $B(R)$. Note also that the function $A(y)$ is
convex. It follows from the fact that in view of the positivity of
the current $\theta$, the measure
$$
4\sum_{j,k}^n\theta'_{j,k}t_j t_k=
\sum_{j,k=1}^n{{\partial^2 A}\over{\partial y_j\partial y_k}}t_jt_k
$$
is positive for any $t_j,t_k\in\R$.

So, the problem of finding a convex function satisfying conditions
(17) is solved for each ball $B(x_0,R)\Subset G$. Denote all such
balls by $B_\alpha$ and the corresponding functions by $A_\alpha$.
Observe that if $B_\alpha\cap B_{\alpha'}\neq\emptyset$, the function
\ben
A_{\alpha',\alpha}:=A_{\alpha'}-A_{\alpha}
\een
is defined on that intersection and satisfies  the equations
$$
{{\partial^2 A_{\alpha',\alpha}}\over{\partial y_j\partial y_k}}=0,
\quad\forall j,k
$$
in the distribution space, so that is a linear function. We also set
$A_{\alpha,\alpha'}=0$ if $B_\alpha\cap B_{\alpha'}=\emptyset$. It is
evident that $A_{\alpha,\alpha'}=-A_{\alpha',\alpha}$ and
$A_{\alpha',\alpha}+A_{\alpha,\alpha"}+A_{\alpha",\alpha'}=0$. Thus
the functions $A_{\alpha,\alpha'}$ form a cocycle in the space of
1-cochains of the covering $\{B_\alpha\}$ with the values in the
bundle of germs of linear functions on $G$. As the domain $G$ is
convex, the corresponding cohomology group equals zero, so that every
cocycle is a coboundary. Therefore there exist linear functions
$h_\alpha(y)$ such that
\ben
h_{\alpha}-h_{\alpha'}=A_{\alpha',\alpha}, \quad\forall\alpha',\alpha.
\een
Set
$$
A(y)=A_{\alpha}(y)+h_{\alpha}(y), \quad y\in B_{\alpha}.
$$
The function $A(y)$ is well defined in the whole domain $G$ due to
(19) and (20). Furthermore, it is convex in $G$ because the
functions $A_\alpha(y)$ are convex in $B_\alpha$ and the functions
$h_\alpha(y)$ are linear. Finally, by the construction of the
functions $A_\alpha(y)$, the equations
$$
{{\partial^2 A_{\alpha}}\over{\partial y_j\partial y_k}}=4\theta'_{j,k},
\quad\forall j,k,
$$
have place in $B_\alpha$, so that the equation (17)
holds in $G$. Therefore, the function $A(y)$ is what we sought. The
theorem is proved.

\medskip
{\it Remark 3.1.} Applying Theorem 3.1 to the current of integration
over an almost periodic divisor $Z$ and remembering the preceding
observations about the density of the divisor, we get the following
statement.

\medskip
{\bf Theorem 3.2.} {\sl Let $Z$ be an almost periodic divisor in
$T_G$.  Then there exists a convex function $A_Z(y)$ in $G$ (the
Jessen function of the divisor $Z$) such that its Riesz measure
$\mu_A$ coincides with the density of the divisor $Z$.}

\medskip
Note that the function $A_Z$ is not uniquely defined by the divisor
$Z$, namely up to a linear term. Besides, it is clear that in the
case of $Z=Z_f$, where $f$ is an almost periodic functionn from
$H(T_G)$, $A_Z=A_f+L$, $L$ being a linear function.

\medskip
{\it Remark 3.2.} An analysis of the proof of Theorem 3.1 shows that
incidentally the following result has been proved.

\medskip
{\bf Proposition 3.1.} {\sl Let $F$ be a closed positive current
with the coefficients $F_{j,k}=\hat F_{j,k}\otimes m_n$, where
$\hat F_{j,k}$ are measures on $G$. Then there exist a convex
function $A(y)$ in $G$ and constants $c_{j,k}\in\R,\:
j,k=1,\ldots,n$, such that}
$$
F=dd^c(u(y)+\sum_{j<k}2c_{j,k}(x_iy_j-y_ix_j)),
$$

\medskip
{\it Remark 3.3.} If $f$ is an almost periodic holomorphic function
in $T_G$, then according to (3) and (2),
$$
\mm(dd^c\log|f|^2)=\tilde a_f^{(1)}=dd^c\tilde A_f^0,
$$
so that the corresponding constants $c_{j,k}=0$. Therefore vanishing
of the constants $c_{j,k}$ is a necessary condition for realizability
of an almost periodic divisor $Z$ as the divisor of an almost
periodic  holomorphic function.  Using of this condition allows us to
show that the periodic divisor
$$
\tilde Z=\sum_{(k_1,k_2)\in\Z^2}(Z-k_1-ik_2)
$$
where $Z=(|Z|,\gamma_Z),\ |Z|=\{(z_1,z_2):\: z_1-iz_2=0\},\
\gamma_Z\equiv 1$, cannot be realized not only as the divisor of an
entire periodic function (as was shown in \cite{keyR6}), but also as the
divisor of an entire almost periodic function. Really, simple
calculations give us
\be
[\tilde Z]=i(dz_1\wedge d\bar z_1+dz_2\wedge d\bar z_2-
dz_1\wedge d\bar z_2+dz_2\wedge d\bar z_1)=\\
dd^c(4(x_1^2+y_1^2)+4(x_1y_2-x_2y_1))
\ee
so that $c_{1,2}\neq 0$.
\footnote
{One can show that for a periodic divisor $Z$, the matrix
$\left(C_{j,k}\right)$ coincides (up to a factor $2\pi$) with
the matrix $\tilde N$ introdused in \cite{keyRRu}.
It was proved there that
vanishing of  the matrix $\tilde N_Z$ is equivalent to realizability
of $Z$ as the divisor of a periodic function. It follows thus from
Remark 3.3 that if a periodic divisor can be realized as the divisor
of no periodic holomorphic function, then it can be realized as the
divisor of no almost periodic holomorphic function, too.}

\section{Divisors with piecewise linear Jessen function}

Here we consider convex piecewise linear  functions in a convex
domain $G$, i.~e. convex functions $A(y),\: y\in G$, such that
\begin{enumerate}
\item the set $\Lambda_A$ of the points of $G$ in whose
neighbourhoods the function $A(y)$ is linear, is dense in $G$;
\item the intersection of $\Lambda_A$ with every domain $G'\Subset G$
consists of a finite number of connected components.
\end{enumerate}

The Jessen function of an almost periodic divisor is by definition convex, however not
every convex function can be the Jessen function of some divisor. A
description was given in \cite{keyR5} for convex piecewise linear  functions
that are  the Jessen functions of almost periodic holomorphic functions. The proof is based on the
connection between  the Jessen function of an almost periodic holomorphic function and the density
of its divisor. As follows from Theorem 3.2, the similar connection
exists between an almost periodic divisor and its Jessen function, too. Therefore the
following theorem is true, and the proof repeats the corresponding
arguments from \cite{keyR5} practically word to word and thus is omitted
here.

\medskip
{\bf Theorem 4.1} {\sl In order that a convex piecewise linear
function $A(y)$ in $G$ be the Jessen function of some almost periodic divisor in $T_G$,
it is nessesary and sufficient that, up to a linear term,
\ben
A(y)=\sum_{j=1}^{\omega}\gamma_j(\< y,\lj\>-h_j)^+,
\een
where $\omega\le\infty,\ \lj\in\R^n,\ |\lj|=1,\
h_j\in\R,\ \gamma_j>0, \ (\cdot )^+=\max (\cdot ; 0)$. Furthermore,
the support $|Z|$ of the divisor is the union of complex hyperplanes
of the form}
\ben
L_{j,p}=\{z:\< z,\lj\>-ih_j-\alpha_{j,p}=0\}
\quad\alpha_{j,p}\in\R.
\een

\medskip
By the help of the results from \cite{keyRaRF} concerning divisors in $\C$,
this theorem can be supplemented by the statement of realizability of
an almost periodic divisor with piecewise linear Jessen function as the divisor of
an almost periodic holomorphic function from a specified class.

\medskip
{\bf Theorem 4.2} {\sl In order that the Jessen function $A_Z(y)$ of
an almost periodic divisor $Z\in T_G$ be piecewise linear,
it is nessesary and sufficient that the divisor $Z$ be the divisor of
a function $f\in H(T_G)$ of the form
\ben
f(z)=\prod_{j=1}^{\omega}f_j(\< z,\lj\>-ih_j),
\een
where $\omega\le\infty,\ \lambda^{(j)}\in\R^n,\ |\lambda^{(j)}|=1,\
h_j\in\R$, and $f_j(w)$ are entire almost periodic functions of $w\in\C$ with real
zeros.}

\medskip
{\it Proof.} If an almost periodic function $F(z)\neq 0$ in a domain $T_{G'},\ G'\subset
G$, the Jessen function $A_F(y)$ is evidently linear in $G'$. Then
the Jessen function of the function $f(z)$ from (23) and thus the
Jessen function of its divisor is linear out of the union of the
planes
\ben
\hat L_j=\{y\in\R^n:\<y,\lj\>-h_j=0\},
\quad j=1,2,\ldots,\omega
\een
Furthermore, given a point $y\in G$, there exists its neighbourhood
which is intersecting with finite many planes (24) only. That proves
the sufficiency part of the statement.

Let now $Z$ be an almost periodic divisor in $T_G$ with piecewise linear Jessen
function $A_Z(y)$. By Theorem 4.1, the function $A_Z(y)$ has the
representation (21) and $|Z|$ is the union of complex hyperplanes
$L_{j,p}$ of the form (22). As was mentioned above, vanishing of the
density of an almost periodic divisor in a domain $G'$ implies the relation $Z\cap
T_{G'}=\emptyset$, so the projections of any hyperplane $L_{j,p}$ to
$\R_{(y)}^n$ with a fixed $j$ have the form $\hat L_j$ from (24)
(with the same $j$).

Note that a transform of the form $z'=B\,z+i\,a$ where $B$ is an
orthogonal matrix with real coefficients and $a\in\R^n$, affects neither
almost periodicity of the divisor, nor piecewise linearity of the
function $A_Z$. Hence one can take $\hat L_j=\{y:\: y_1=0\}$ and,
respectively, $L_{j,p}=\{z:\: z_1=\alpha_{j,p}\}$, where
$\alpha_{j,p}$ are the same as before. One can also take $0\in\hat
L_j$ belongs to no other plane $\hat L_{j'},\: j'\neq j$.

Let $\gamma _{j,p}$ be the density of the divisor in a regular point
of $L_{j,p}$. It remains evidently the same for all the points that
are not the points of intersection of  $L_{j,p}$ with
$L_{j',p},\: j'\neq j$. Choose $\e>0$ in a way that the
interval $I=\{y\in\R^n:\: |y_1|<\e,\ y_2=\ldots =y_n=0\}$ is
contained in $G$ and intersects no plane $L_{j',p},\: j'\neq j$.
Consider the divisor $Z_j$ in the strip $\Pi_\e =\{w\in\C:\:
|{\rm Im}\, w|<\e\}$ formed by the points $\alpha_{j,p}$  with
the multiplicities $\gamma_{j,p}$. We claim that this divisor is
almost periodic. To prove it, notice that in the one-dimensional situation the
integration current over any divisor $Z$ is the measure
$\delta_Z$ concentrated on the support of the divisor and such that
the measure of each point equals its multiplicity. Therefore it
suffices to prove that the function
$$
(\delta_{Z_j}(w),\psi(w-t))=\sum_p\gamma_{j,p}\psi(\alpha_{j,p}-t)
$$
is almost periodic on $\R$ for any function $\psi\in{\cal D}(\Pi_\e)$ with
the support in a small enough neighbourhood of the origin.

By almost periodicity of the divisor $Z$, the function
\ben
(\tr[Z],\varphi(z-t))=\sum_{j,p}\gamma_{j,p}\int_{L_{j,p}}
\varphi(z-t)dV_{j,p}(z),\quad t\in\R^n,
\een
where $dV_{j,p}$ is the elememt of $(2n-2)$-volume on $L_{j,p}$, is
an almost periodic function on $\R^n$ for any $\varphi\in{\cal D}(T_G)$. 
Choose a function $\varphi(z)$ of the form
$\varphi(z)=\varphi_1(z_1)\,\varphi_2('z)$, where $'z=(z_2,\ldots,z_n)$,
such that $\supp\varphi_1\subset\Pi_\e$ and $\supp\varphi_2
\subset \{'z\in\C^{n-1}:\: |'z|<\e_1\}$ with $\e_1$ small
enough that the projection of the set $\supp\varphi_1\times
\supp\varphi_2$ to $\R_{(y)}^n$ is out of the planes
$\hat L_{j'}$ with $ j'\neq j$. Then the function in (25) takes the
form
$$
F_{\varphi}(t_1,'t)=\sum_{p=1}^{\infty}\gamma_{j,p}\int_{L_{j,p}}
\varphi_1(z_1-t_1)\varphi_2('z-'t)dV_{j,p}('z),
$$
where $'t=(t_2,\ldots,t_n)$. Since this function is almost periodic in $\R^n$,
the function
$$
F_{\varphi}(t_1,'0)=
C\sum_{p=1}^{\infty}\gamma_{j,p}\varphi_1(\alpha_{j,p}-t_1),
$$
where
$$
C=\int_{\C_{n-1}}\varphi_2('z)dV_{2n-2}('z).
$$
is an almost periodic function on $\R$. As the function $\varphi_1\in{\cal
D}(\Pi_\e)$ has been chosen arbitrarily, almost periodicity
of the divisor $Z_j$ is proved.

Observe now that since ${\rm Card}\,\{j:\: G'\cap\hat L_{j}\neq
\emptyset\}<\infty$ for every $G'\Subset G$, one can choose a
sequence of domains $G_k\Subset G_{k+1},\ k=1,2,\ldots,\ G=\cup G_k$,
such that $\bar G_k\cap \hat L_{k}=\emptyset$. Replacing if necessary
$\lambda^{(k)}$ and $h_k$ with $-\lambda^{(k)}$ and $-h_k$ one may
take
$$
\bar G_k\subset\{y\in\R^m:\<y,\lambda^{(k)}\>-h_k>0\}
$$
and thus, for some $0<r_k<R_k<\infty$ and any $l\le k$,
\ben
\bar G_l\subset\bar G_k\subset\{y\in\R^m:r_k<\<y,\lambda^{(k)}\>-h_k<R_k\}
\een

Now we can use Theorem 2 from \cite{keyRaRF} resulting 
that every almost periodic 
divisor in $\C$ with the support in $\R$ is realizable 
as the divisor of an entire almost periodic function. 
So, let $\tilde f_j(w)$ be an entire almost periodic function in $\C$ with
$Z_{\tilde f_j}=Z_j$. Since $\tilde f_j(w)$ does not vanish out of
the real axis, one has in $\{w:\: {\rm Im}\,w>0\}$ the form (see for example
\cite{keyLt})
$$
\tilde f_j(w)=\exp\{ic_jw+g_j(w)\},
$$
where $c_j\in\R$, and $g_j(w)$ are almost periodic holomorphic
functions in the upper half-plane.  By the approximation theorem for
almost periodic functions, there exist finite exponential sums
$q_j(s)$ such that
\ben
\sup\{|g_j(w)-q_j(w)|:r_j<{\rm Im}w<R_j\}\le j^{-2}.
\een

Now set
$$
f_j(w)=\tilde f_j(w)\exp\{-ic_jw-q_j(w)\}.
$$
It follows from (27) and (26) that for any fixed $k$, the
inequality
$$
\sum_{j=k}^{\omega}|\log f_j(\<\lambda^{(j)},z\>-ih_j)|\le
\sum_{j=k}^{\omega}j^{-2}
$$
takes place in the domain $G_k$. Therefore the product (23)
converges uniformly on every domain $T_{G'}$ with $G'\Subset G$, and
is an almost periodic function in $T_G$ because each function
$f_j(\langle\lambda^{(j)},
z\rangle -i\,h_j)$ is almost periodic in $\C^n$.

The rest follows from the fact that the zero set of the
function $f_j(\langle\lambda^{(j)},z\rangle -i\,h_j)$ coincides with
$\cup_{p=1}^\infty L_{j,p}$, and if a point $z\in L_{j,p}$ is a
regular point of the zero set, the multiplicity of this point as a
zero of the function equals $\gamma_{j,p}$. Therefore the divisor
$Z_f$ coincides with the given divisor $Z$. The proof is complete.


\begin{thebibliography}{A}

\bibitem{keyF}
S.Yu.\,Favorov, {\it Estimates for asymptotic measures and Jessen
functions for almost periodic functions,}  Dopov. Nats. Akad. Nauk
Ukraini, {\bf 10}, (1996), 27-30  (russian). 


\bibitem{keyGrK}
F.A.\,Griffiths and J.\,King, {\it Nevanlinna theory and holomorphic mappings
between algebraic varieties,} Acta Math. {\bf 130} (1973), 145-220.

\bibitem{keyH}
R.\,Harvey, {\it Holomorphic Chains and their Boundaries,}
Proceedings of Symposia in Pure Mathematics, {\bf 30}, part 1,
(1977), 307-382.

\bibitem{keyJ}
B.\,Jessen, {\it Ueber die Nullstellen einer analitischen fastperiodischen
Functions,} Eine Verallagemeinerung der Jensenschen Formel.
Math. Ann., {\bf 108}, (1933), 485-516.

\bibitem{keyKL}
M.G.\,Krejn and B.Ja.\,Levin, {\it On almost periodic functions of
exponential type}, Dokl. AN SSSR, {\bf 64} (1949), no.3,
285-287.

\bibitem{keyLa}
N.S.\,Landkof, {\it Foundations of Modern Potential Theory,}
Springer-Verlag, Berlin, 1972.

\bibitem{keyLe}
B.Ja.\,Levin, {\it Distribution of Zeros of Entire Functions.} Transl.
of Math. Monograph, Vol.5, AMS Providence, No.1 (1964), MR 19, 403,
494 p.

\bibitem{keyLt}
B.M.\,Levitan, {\it Almost Periodic Functions}, FML, Moskow,
1953, 396 p.

\bibitem{keyLG}
P.\,Lelong, L.\,Gruman, {\it Functions of Several Complex Variables,}
Springer-Verlag. Berlin-Heidelberg (1986), 270p.

\bibitem{keyRa1}
A.Yu.\,Rashkovskii, {\it Currents associated to holomorphic almost periodic
mappings,}
  Mat. Fizika, Analiz i Geometria, 2 (1995), No. 2, 150-169.

\bibitem{keyRa2}
A.Yu.\,Rashkovskii, {\it Monge-Amp\`ere type currents associated to
holomorphic almost periodic mappings,} C. R. Acad. Sci. Paris, S\'erie I
{\bf 321} (1995), 1553-1558.

\bibitem{keyRaRF}
A.Yu.\,Rashkovskii, L.I.\,Ronkin and S.Yu.\,Favorov,
{\it On Almost Periodic Sets in the Complex Plane,}
Dopov. Nats. Akad. Nauk Ukra\"i ni, (russian), to appear.

\bibitem{keyR4}
L.I.\,Ronkin, {\it Jessen's theorem for holomorphic almost periodic
functions in tube domains,} Sibirsk. Mat. Zh. {\bf 28} (1987), 199-204.

\bibitem{keyR2}
L.I.\,Ronkin, {\it Jessen's theorem for holomorphic almost periodic mappings,}
 Ukrainsk. Mat. Zh. {\bf 42} (1990), 1094-1107 (russian)

\bibitem{keyR3} L.I.\,Ronkin {\it Functions of Completely Regular Growth}.
Mathematics and its Applications, Soviet Ser., vol.81
Kluver Academic Publ., Dordrecht Boston, 1990, 392 p.

\bibitem{keyR5}
L.I.\,Ronkin, {\it On a certain class of holomorphic almost periodic
functions,}
 Sibirskii Mat. Zh. {\bf 33} (1992), 135-141, (russian).

\bibitem{keyR6}
L.I.\,Ronkin, {\it Holomorphic periodic functions and periodic divisors,}
Math. Physics, Anal. and Geom. {\bf 2} (1995),108-122 (russian).

\bibitem{keyR1}
L.I.\,Ronkin, {\it Almost periodic distributions and divizors
in tube domains}, Zap. Nauchn. Sem. POMI {\bf 247} (1997).

\bibitem{keyRRu}
L.I.\,Ronkin and A.M.\,Russakovskii, {\it Entire Periodic Functions 
with Plane Zeros,} Complex Variables, {\bf 30}, (1996), 323-338.

\bibitem{keyT1}
H.\,Tornehave, {\it Systems of zeros of holomorphic almost periodic
 functions}, Kobenhavns Universitet Matematisk Institut, Preprint No. 30,
(1988), 52 p.

\bibitem{keyT2}
H.\,Tornehave, {\it On the zeros of entire almost periodic function},
The Harald Bohr Centenary (Copenhagen 1987). Math. Fys. Medd. Danske,
{\bf 42} (1989), no. 3, 125-142.

\end{thebibliography}
\end{document}